% Commenced on 5 May 2015,

%%==== PLAIN TeX with AMS fonts & symbols ====%%

%\magnification=1200
\nopagenumbers
\parindent= 15pt
\baselineskip=14pt

%-----Paper Format---------
\hsize=15cm
\vsize=23cm
\hoffset=0.55cm
\voffset=0.3cm

%-----AMS Symbols-------
\input amssym.def
\input amssym.tex

%-----Special Fonts-----
 at6.5pt
\font\srm=cmr8

\font\csc=cmcsc10
\font\title=cmr12 at 14pt
%\font\stitle=cmsl12 at 14pt

\font\teneusm=eusm10    
\font\seveneusm=eusm7  
\font\fiveeusm=eusm5    
\newfam\eusmfam

\textfont\eusmfam=\teneusm 
\scriptfont\eusmfam=\seveneusm
\scriptscriptfont\eusmfam=\fiveeusm

%----Small Macros----

\def\r#1{{\rm #1}}

\def\B#1{{\Bbb #1}}

%----Headers-------
\def\rightheadline{\hfil{\srm Cornacchia's algorithm}
\hfil\tenrm\folio}
\def\leftheadline{\tenrm\folio\hfil{\srm Y. Motohashi}\hfil}
\def\emptyheadline{}
\headline{\ifnum\pageno=1 \emptyheadline\else
\ifodd\pageno \rightheadline \else \leftheadline\fi\fi}

%-----Footers------
\def\firstpage{\hss{\vbox to 2cm
{\vfil\hbox{\rm\folio}}}\hss}
\def\emptyfootline{\hfil}
\footline{\ifnum\pageno=1\firstpage\else
\emptyfootline\fi}

%----Text---------
\centerline{\title On Cornacchia's Algorithm}
\bigskip
\centerline{\csc Yoichi Motohashi}
\vskip 0.7cm 
\noindent
{\bf Abstract:}  We give an endorsement for Cornacchia's
famous algorithm. Thus we do not claim anything new but an
approach which is supposed to be simpler
than those of previous works written with the same aim. 
All variables and constants are integers.
\bigskip
\noindent
{\bf 1. Introduction}
\smallskip
We consider the Diophantine equation
$$
x^2+dy^2=m,\eqno(1.1)
$$
where
$$
d\ge2,\quad m\ge2,\quad \r{gcd}\{d,m\}=1.\eqno(1.2)
$$
We observe that provided (1.1) admits a proper solution $\{u,v\}$, 
namely
$$
u^2+dv^2=m,\quad \r{gcd}\{u,v\}=1,\eqno(1.3)
$$
it holds that
$$
\r{gcd}\{uv,m\}=1,\quad uv\ne0,\quad m\ge d+1;\eqno(1.4)
$$
and that $-d$ is a quadratic residue $\bmod\, m$, and there exists 
a $w$, $m/2\le w<m$, such that
$$
w^2\equiv-d\bmod m.\eqno(1.5)
$$
In fact one may put $w\equiv uv^{-1}\bmod m,\; 
vv^{-1}\equiv1\bmod m$.
\smallskip
We apply the antenaresis to the pair $\{w,m\}$ and denote by 
$\{t_j\}$ the residues thus arising. Also we expand $w/m$ into
a regular continued fraction and get the convergents
$\{C_j/D_j\}$. An explanation of these notions is
to be given in the next section.
\smallskip
With this, Cornacchia (1908, pp.\ 60--66) essentially stated the
following:
\medskip
\noindent
{\bf Theorem.} 
\smallskip
\noindent
{\it Let $\nu$ be such that
$$
t_{\nu+1}^2\le m<t_\nu^2.\eqno(1.6)
$$ 
Then $(1.3)$ implies that
$$
|u|=t_{\nu+1},\quad |v|=D_\nu.\eqno(1.7)
$$
}
\noindent
{\bf Algorithm.}
\smallskip
\noindent
{\it
Solely on the assumption $(1.5)$, that is, without
any prior knowledge of $(1.3)$, find $\nu$ that satisfies $(1.6)$.
If it holds that $t_{\nu+1}^2+dD_\nu^2=m$, then this is a
proper solution of $(1.1)$ corresponding to $w$. 
Otherwise $(1.1)$ does not admit any proper solution  
corresponding to $w$.
}
\medskip
Obviously, if $(1.5)$ is empty, then there is no need to make any 
further quest as far as proper solutions are concerned; but
see the second example in the last section.
We shall give, in the third section, 
a proof of the theorem and algorithm 
using a basic observation on the nature of
finite continued fractions, which is in fact contained in the best approximation
theorem of Lagrange (1798, pp.\ 55--57) but 
can be proved quickly with
an idea of Legendre (1798, p.\ 27). Also,
we partly follow Basilla (2004). The case $d=1$ will be treated in the
fourth section.
\medskip
\noindent
{\it Acknowledgements\/}: We are indebted to Professor 
Francesco Pappalardi who
kindly provided us with a copy of Cornacchia's work which would have
been hard for us to acquire otherwise. Also, we are grateful to 
Nigel Watt for his helpful comments on a draft version of the present note.
\bigskip
\noindent
{\bf 2. Finite continued fractions}
\smallskip
Let $b\ge2$ and expand any fraction $a/b$ into the
continued fraction
$$
q_0+{1\over q_1}{\atop +}{1\over q_2}
{\atop+\cdots}{\atop+}{1\over q_k},\eqno(2.1)
$$
which means that the antenaresis applied to $\{a,b\}$ yields the
sequence of identities
$$
r_{j-1}=q_jr_j+r_{j+1},\quad 0\le r_{j+1}<r_j,\quad0\le j\le k,
\eqno(2.2)
$$
with the convention $r_{-1}=a,\, r_0=b,\, r_{k+1}=0$; in particular
$r_k=\r{gcd}\{a,b\}$. One may put this,
in the matrix multiplication format, as
$$
(r_{j-1},r_j)=(r_j,r_{j+1})\left(\matrix{q_j&1\cr1&0}\right),
\quad0\le j\le k.\eqno(2.3)
$$
Writing, for $-1\le j\le k$,
$$
\eqalignno{
&\qquad\left(\matrix{A_j & A_{j-1}\cr B_j & B_{j-1}}\right)
=\left(\matrix{q_0&1\cr 1&0}\right)
\left(\matrix{q_1& 1\cr 1&\hfil 0}\right)\cdots
\left(\matrix{q_j&1\cr 1&0}\right),&(2.4)\cr
& A_{-2}=0,\; B_{-2}=1;
\quad A_{-1}=1,\; B_{-1}=0;\quad A_0=q_0,\; B_0=1,&(2.5)
}
$$
where an empty product is the unit matrix, we have, by induction,
$$
\eqalignno{
&\;{A_j\over B_j}=q_0+
{1\over q_1}{\atop +}{1\over q_2}
{\atop+}{\atop\cdots}{\atop+}{1\over q_j},\quad
0\le j\le k,&(2.6)\cr
&A_jB_{j-1}-A_{j-1}B_j=(-1)^{j-1},\quad -1\le j\le k.&(2.7)
}
$$
Thus, we have
$$
\eqalignno{
(r_{-1},r_0)&\,=(r_j,r_{j+1})\left(\matrix{q_j&1\cr1&0}\right)
\left(\matrix{q_{j-1}&1\cr1&0}\right)\cdots 
\left(\matrix{q_0&1\cr1&0}\right)\cr
&\,=(r_j,r_{j+1})\left(\matrix{A_j&B_j\cr A_{j-1}& B_{j-1}}
\right),&(2.8)
}
$$
in which a use is made of the fact that 
$\left({q_j\atop1}{1\atop0}\right)$ 
are symmetric.
Hence, we get
$$
(r_j,r_{j+1})=(-1)^{j-1}(a,b)
\left(\matrix{\hfill B_{j-1}& -B_j\cr -A_{j-1}& \hfill A_j}\right).
\eqno(2.9)
$$
\medskip
\noindent
{\bf Lemma.}
\smallskip
\noindent
{\it If it holds, with a particular $\lambda$,
$-1\le \lambda\le k$, that
$$
|aQ-bP|<r_\lambda,\quad Q\ne0,\eqno(2.10)
$$
then we have
$$
B_\lambda\le|Q|.\eqno(2.11)
$$
}
\par
\noindent
{\it Proof\/}. Obviously we may skip the cases $\lambda=-1,0$;
we assume $1\le \lambda\le k$. Following Legendre {\it loc.cit.\/}, 
we introduce the transformation of variables
$$
P=MA_\lambda-NA_{\lambda-1},\;
Q=MB_\lambda-NB_{\lambda-1}; \eqno(2.12)
$$
namely 
$$
M=(-1)^\lambda(QA_{\lambda-1}-PB_{\lambda-1}),\;
N=(-1)^\lambda(QA_\lambda-PB_\lambda).\eqno(2.13)
$$ 
We have, via $(2.9)$,
$$
\eqalignno{
aQ-bP&\,=M(aB_\lambda-bA_\lambda)
-N(aB_{\lambda-1}-bA_{\lambda-1})\cr
&\,=(-1)^\lambda(Mr_{\lambda+1}+Nr_\lambda).&(2.14)
}
$$
If $MN>0$, then $|aQ-bP|=|M|r_{\lambda+1}+|N|r_\lambda
\ge r_\lambda$, which is rejected by the assumption $(2.10)$. If
$M=0$, then $(2.13)$ implies $Q=\tau B_{\lambda-1},\, 
P=\tau A_{\lambda-1}$ with a $\tau\ne0$, since $\r{gcd}\{B_j,
A_j\}=1$ by $(2.7)$; hence, via $(2.9)$, we get
$|aQ-bP|=|\tau|r_\lambda\ge r_\lambda$, 
which contradicts $(2.10)$.
Therefore, we may suppose that
$MN\le0,\, M\ne0$, and find that $(2.12)$ implies
$|Q|=|M|B_\lambda+|N|B_{\lambda-1}\ge B_\lambda$. 
We end the proof.
\bigskip
\noindent
{\bf 3. Proof of Cornacchia's theorem and algorithm}
\smallskip
We specialize the discussion of the 
previous section by setting $a=w$, $b=m$,
$r_j=t_j$, $A_j=C_j$, $B_j=D_j$. By definition,
$t_0=m,\, t_1=w,\, t_k=1$; thus there exists a unique $\nu$ that
satisfies $(1.6)$. On the other hand, $(2.9)$ gives that
$t_{j+1}=(-1)^j(wD_j-mC_j)$ and
$$
t_{j+1}^2+dD_j^2\equiv (w^2+d)D_j^2\equiv0\bmod m,
\quad -1\le j\le k.\eqno(3.1)
$$
Also, we have $u=vw-\ell m$ with an $\ell$; and
by the lemma with $P=\ell,\, Q=v$ we see that
$$
\hbox{for any $\lambda$ such that 
$|u|< t_\lambda$,  $-1\le \lambda\le k$,}\atop
{\hbox{we have $D_\lambda\le|v|$, that is,  $dD_\lambda^2<m$.}}
\eqno(3.2)
$$
In particular, since $|u|<\sqrt{m}<t_\nu$, we have
$t_{\nu+1}^2+dD_\nu^2<2m$DHence by the congruence relation
$(3.1)_{j=\nu}$, we find that
$$
t_{\nu+1}^2+dD_\nu^2=m.\eqno(3.3)
$$
If $|u|>t_{\nu+1}$, then $m=u^2+dv^2>t_{\nu+1}^2
+dD_\nu^2=m$, which is impossible. 
If $|u|<t_{\nu+1}$, then in much the same way as above
we get $t_{\nu+2}^2+dD^2_{\nu+1}=m$; and
$$
\eqalignno{
m^2&\,=\big|t_{\nu+1}+iD_\nu\sqrt{d}\big|^2
\big|t_{\nu+2}+iD_{\nu+1}\sqrt{d}\big|^2\cr
&\,=(t_{\nu+1}t_{\nu+2}-dD_\nu D_{\nu+1})^2
+d(t_{\nu+1}D_{\nu+1}+t_{\nu+2}D_\nu )^2,&(3.4)
}
$$
which is impossible, since $d\ge2$ and $m=t_{\nu+1}D_{\nu+1}
+t_{\nu+2}D_\nu$ by $(2.8)_{j=\nu+1}$. Hence,
$(1.7)$ is verified. We end the proof of the theorem.
\smallskip
As to the validity of the algorithm,
 it suffices to show that the identity $(3.3)$, if it
holds on its own, implies that $\r{gcd}\{t_{\nu+1}, D_\nu\}=1$. 
To see this, we put $w^2+d=hm$, and get,
via $(2.9)_{j=\nu}$,
$$
t_{\nu+1}^2+dD_\nu^2=m\big(hD_\nu^2-2wC_\nu D_\nu
+mC_\nu^2\big).\eqno(3.5)
$$
Namely, $(hD_\nu-2wC_\nu)D_\nu+mC_\nu^2=1$. Hence
$\r{gcd}\{D_\nu,m\}=1$, and
$\r{gcd}\{t_{\nu+1},D_\nu\}=1$. This ends the 
endorsement of Cornacchia's algorithm.
\bigskip
\noindent
{\bf 4. Remarks and examples}
\smallskip
\noindent
{\it Remark\/} 1: When $d=1$ we follow the argument 
of Hermite (1848). 
Thus we assume {\it only\/} that there 
exists a $w$, $m/2\le w<m$, such that
$$
w^2\equiv-1\bmod m;\eqno(4.1)
$$
and we adopt the specialization at the beginning of the
last section. Then we choose $\mu$ to satisfy
$$
D_\mu^2\le m<D_{\mu+1}^2.\eqno(4.2)
$$
This is possible because $\{D_j\}$ increase monotonically from $1$ to
$m$.  We note that, for $1\le j\le k$,
$$
\left|{w\over m}-{C_{j-1}\over D_{j-1}}\right|
\le{1\over D_{j-1} D_j}
\quad\hbox{or\quad $t_j\le m/D_j$},\eqno(4.3)
$$
since by the construction $w/m$ is 
between $C_{j-1}/D_{j-1}$ and $C_j/D_j$, and 
we have $(2.7)$ and $(2.9)$. Hence
we see that $t_{\mu+1}<\sqrt{m}$,
and $t_{\mu+1}^2+D_\mu^2<2m$. The congruence relation
$(3.1)$ holds when $d=1$, as well; and it implies
$$
t_{\mu+1}^2+D_\mu^2=m.\eqno(4.4)
$$
As in the case $d\ge2$, this is a proper solution of $(1.1)_{d=1}$.
\smallskip
Then we assume that we have a
proper solution $(1.3)_{d=1}$. We
observe that
$$
\eqalign{
&vD_\mu-ut_{\mu+1}\equiv vD_\mu-(-1)^\mu vw^2 D_\mu
\equiv (1+(-1)^\mu)vD_\mu\bmod m,\cr
&uD_\mu-vt_{\mu+1}\equiv uD_\mu
-(-1)^\mu vwD_\mu\equiv (1-(-1)^\mu)uD_\mu\bmod m,
}\eqno(4.5)
$$
one of which is congruent to $0\bmod m$.
Let us assume $vD_\mu-ut_{\mu+1}\equiv0\bmod m$.
Then we note that
$$
(vD_\mu-ut_{\mu+1})^2+(uD_\mu+vt_{\mu+1})^2=
|u+iv|^2|t_{\mu+1}+iD_\mu|^2=m^2.\eqno(4.6)
$$
This implies that either $(vD_\mu-ut_{\mu+1})^2=0$ or 
$=m^2$; and if the
latter holds, then $(uD_\mu+vt_{\mu+1})^2=0$.
Since $\r{gcd}\{u,v\}=1$ and $\r{gcd}\{t_{\mu+1},D_\mu\}=1$,
we conclude that
$$
\hbox{either\quad $|u|=D_\mu,\, |v|=t_{\mu+1}$\quad or \quad
$|u|=t_{\mu+1},\, |v|= D_\mu$.}\eqno(4.7)
$$
With the remaining case, i.e.,
$uD_\mu-vt_{\mu+1}\equiv0\bmod m$, we use
instead $|u-iv|^2|t_{\mu+1}+iD_\mu|^2=m^2$, 
getting $(4.7)$ again.
\smallskip
It should be stressed that in the 
case $d=1$ we do not need to have $(1.3)$; it suffices to
have $(1.5)_{d=1}$ or $(4.1)$. However, 
we should first choose $D_\mu$ instead of $t_{\mu+1}$.
\medskip
\noindent
{\it Remark 2\/}: It is worth remarking that 
Smith (1855) showed that when dealing with $m$ a prime
$\equiv 1\bmod 4$ the condition $(1.5)_{d=1}$ is not needed to
be assumed as far as one is concerned with only the existence of
the representation $(1.3)_{d=1}$.
He exploited the fact that there exists a fraction $m/h$, $2\le h<m/2$,
whose continued fraction expansion is palindromic. 
\medskip
\noindent
{\it Remark 3\/}: An effective way to adopt prior to 
any use of Cornacchia's algorithm is to restrict oneself 
to the cases of $m$ being square-free and
prime powers. Then in the general case
the identity $|x+iy\sqrt{d}|^2=x^2+dy^2$ is to be exploited.
It is important not to restrict oneself to prime powers only, 
since the product
of primes, none of which has the representation $(1.3)$,
may admit the representation. An example is given below.
\medskip
\noindent
{\it Example 1\/}: We consider the situation $d=5$ and $m=
435629$. Since $2^{m-1}\not\equiv1\bmod m$, we see that
$m$ is a composite number. The $\rho$ method of Pollard, for instance,
gives the decomposition $m=p_1p_2$, with primes $p_1=367,\,
p_2=1187$. 
Since $p_1, p_2\equiv 7\bmod 20$, they are not
expressible by the quadratic form $x^2+5y^2$, but the product
$m=p_1p_2$ does admit a proper representation by the form,
according to a well-known criterion which can be traced back to
Fermat (1658: 1894, p.\ 405).  Let us confirm this by
Cornacchia's algorithm. Following Lagrange (1768, p.\ 500),
we put $w_l\equiv \pm(-5)^{(p_l+1)/4}\bmod p_l$, 
$l=1,2$; namely
$w_1\equiv\pm27\bmod p_1,\, w_2\equiv\pm282\bmod p_2$.
Then $w_l^2\equiv-5\bmod p_l$.
Checking that
$$
{367\over 1187}=0+{1\over 3}{\atop+}{1\over4}{\atop+}
{1\over3}{\atop+}{1\over1}{\atop+}{1\over2}{\atop+}{1\over1}
{\atop+}{1\over5}\;\Rightarrow\;207\cdot 367-64\cdot1187=1,
\eqno(4.7)
$$
we have that
$w\equiv207\cdot 367\cdot w_2-64\cdot1187
\cdot w_1\bmod m$ satisfies $w^2\equiv-5\bmod m$, or
$w=231183,\,386057$.
Then the algorithm yields the solutions
$t_{\nu+1}^2+5D_\nu^2=m$ with
$$
\eqalign{
&w=231183:\quad t_\nu=1385,\;t_{\nu+1}=228,\; 
D_\nu=277,\cr
&w=386057:\quad t_\nu=1450,\;t_{\nu+1}=123,
\; D_\nu=290.
}\eqno(4.8)
$$
Further, as to the equation $x^2+5y^2=p_1^2$, we use the
Sch\"onemann--Hensel lifting, and from $27\bmod p_1$ obtain 
$w_*=109760$, $w_*^2\equiv-5\bmod p_1^2$. The antenaresis 
applied to the pair $\{w_*, p_1^2\}$ yields $362^2
+5\cdot27^2=p_1^2$DWith this, we exploit
$(362+27\sqrt{5}i)(228-277\sqrt{5}i)$, and find
the following proper representation of $p_1^3p_2$:
$$
119931^2+5\cdot94118^2=58674434381.\eqno(4.9)
$$ 
The other combination $(362+27\sqrt{5}i)(228+277\sqrt{5}i)$
does not lead to a proper representation, a phenomenon
which can be explained by a use of the theory of ideals in
$\B{Q}(\sqrt{-5})$.
\medskip
\noindent
{\it Example 2\/}: Consider $x^2+7y^2=4p$ where
$p=9241$ is a prime. This equation does not have any
proper solution, since if $x,y$ are both odd, then $x^2+7y^2
\equiv 0\bmod 8$, and if there is a solution, then $x,y$ are both
even. Thus, Cornacchia's algorithm should not work with this situation,
even though we have $(1.5)$, that is, $w^2\equiv-7\bmod 4p$,
$w=24899$. Applying the antenaresis to the pair $\{w,4p\}$,
we find that $t_{\nu+1}=52,\, D_\nu=144$, and
$t_{\nu+1}^2+7D_\nu^2\ne4p$ indeed. Hence, one should follow the
prescription indicated in the third remark above. We have
$6417^2\equiv-7\bmod p$, and the algorithm gives
$t_{\nu+1}=13,\,D_\nu=36$, and 
$13^2+7\cdot 36^2=p$, which yields the
improper solution $26^2+7\cdot 72^2=4p$.
By the way, the quadratic
residue $6417\bmod p$ is readily obtained 
with the algorithm of Tonelli (1891).
\vskip 1cm
\centerline{\bf References}
\medskip
\noindent
\item{[1]}
J.M. Basilla (2004): On the solution of $x^2+dy^2=m$.
Proc.\ Japan Acad., {\bf 80}(A), 40--41.
\medskip
\noindent
\item{[2]} G. Cornacchia (1908): Su di un metodo per la risoluzione 
in numeri interi dell'equazione $\sum_{h=0}^n C_hx^{n-h}y^h =P$. 
Giornale di matematiche di Battaglini, {\bf 46}, 33--90.
\medskip
\noindent
\item{[3]} P. de Fermat (1658): XCVI, a Kenelm Digby. 
In: {\OE}uvres, tome deuxi\`eme, pp.\ 402--408. 
Gauthier--Villars, Paris 1894.
\medskip
\noindent
\item{[4]} C. Hermite (1848): Note.
J.\ math.\ pures et appliq., {\bf 13}, 15.
 \medskip
\noindent
\item{[5]} J. Lagrange (1769): Sur la solution des probl\`emes 
ind\'etermin\'es du second degr\'e.  In: {\OE}uvres 2, pp.\ 377--535.
Gauthier-Villars, Paris 1868.
\medskip
\noindent
\item{[6]} --- (1798): Additions aux \'el\'emants d'alg\`ebre
d'Euler. Analyse ind\'etermin\'ee.  In: {\OE}uvres 7, pp.\ 5--180.
Gauthier-Villars, Paris 1877.
\medskip
\noindent
\item{[7]}
A.-M. Legendre (1798 (An VI)): Essai sur la th\'eorie 
des nombres. Duprat, Paris.
\medskip
\noindent
\item{[8]} H.J.S. Smith (1855): De compositione 
numerorum primorum formae $4\lambda+1$ ex duobus quadratis. 
J. reine angew.\ Math., {\bf 50}, 91--92. 
\medskip
\noindent
\item{[9]} A. Tonelli (1891): Bemerkung \"uber die Aufl\"osung 
quadratischer Congruenzen. 
Nach\-richten\ K\"onigl.\ Gesell.\ Wiss.\ Georg-Augusts-Univ.\ G\"ottingen, 
344--346.\bigskip
\hfill Email: ymoto@math.cst.nihon-u.ac.jp

\bye